\documentclass[12pt,a4paper]{article}
\usepackage{amsfonts}
\usepackage{amsmath}
\usepackage{mathrsfs}
\usepackage{graphicx}
\usepackage{amsmath,amsthm,amssymb,amscd}

\newtheorem{Thm}{Theorem}[section]
\newtheorem{Lem}[Thm]{Lemma}
\newtheorem{Prop}[Thm]{Proposition}

\newtheorem{Cor}[Thm]{Corollary}
\theoremstyle{remark}
\newtheorem{Rem} [Thm]{Remark}

\theoremstyle{claim}

\DeclareMathOperator{\Diff}{Diff}

\def\m{{\cal M}}

\begin{document}

\begin{center}
{\Large \bf  Ergodic properties of invariant measures for $C^{1+\alpha}$ non-uniformly hyperbolic systems
}\\
\smallskip
\end{center}

\bigskip
\begin{center}
Chao Liang$^*$, Wenxiang Sun$^\dag$ and Xueting Tian$^\ddag$

\smallskip
$^*$ Applied Mathematical Department, The Central University of Finance
and Economics, Beijing 100081, China\\
(chaol@cufe.edu.cn)\\

$^\dag$ LMAM, School of Mathematical Sciences, Peking University, Beijing 100871, China\\
(sunwx@math.pku.edu.cn)

$^\ddag$ Academy of Mathematics and Systems Science, Chinese Academy of Sciences, Beijing 100190, China\\$\&$  School of Mathematical Sciences, Peking University, Beijing 100871, China\\
(tianxt@amss.ac.cn $\&$ txt@pku.edu.cn)
\end{center}
\bigskip

\footnotetext {$^*$Liang is supported by NNSFC(\# 10901167 and \#
10671006)} \footnotetext {$^\dag$Sun is supported by NNSFC(\#
10231020) and Doctoral Education Foundation of China}
\footnotetext {$^\ddag$Tian is the corresponding author}
\footnotetext {AMS Review: 37C40; 37D25; 37H15; 37A35}

\smallskip
\begin{abstract}
For an ergodic hyperbolic measure $\omega$ of a $C^{1+\alpha}$
diffeomorphism, there is a $\omega$ full-measured set
$\tilde{\Lambda}$ such that every nonempty, compact and connected
subset $V$ of $\m_{inv}(\tilde{\Lambda})$ coincides with the
accumulating  set of  time averages of Dirac measures supported at
{\textit{one orbit}}, where $\m_{inv}(\tilde{\Lambda})$ denotes the
space of invariant measures supported on $\tilde{\Lambda}$. Such
state points corresponding to a fixed $V$ are dense in the support
$supp(\omega)$. Moreover $\m_{inv}(\tilde{\Lambda})$ can be
accumulated by  time averages of Dirac measures supported at
{\textit{one orbit}}, and such state points form a residual subset
of $supp(\omega)$. These extend  results of Sigmund \cite{Sig} from
uniformly hyperbolic case to non-uniformly hyperbolic case. As a
corollary,  irregular points form a residual set of $supp(\omega)$.

\end{abstract}

\section{Introduction}

Sigmund \cite{Sig} in 1970 invented  two approximation properties
for $C^1$ uniformly hyperbolic diffeomorphisms: one is that
invariant measures can be approximated by periodic measures, the
other is that every nonempty, compact and connected subset of the
space of invariant measures coincides with  the  accumulating set of
time averages of Dirac measures supported at one orbit and such
orbits are dense. The first approximation property had realized
among $C^{1+\alpha}$ non-uniformly hyperbolic diffeomorphisms in
2003, when Hirayama \cite{Hir} proved that periodic measures are
dense in the set of invariant measures supported on a full measure
set with respect to a hyperbolic $mixing$ measure. In 2009, Liang,
Liu and Sun \cite{LLS} replaced the assumption of hyperbolic
$mixing$ measure by a more natural and weaker assumption of
hyperbolic $ergodic$ measure and generalized Hirayama's result. The
proofs in \cite{Hir,LLS} are both based on Katok's closing and
shadowing lemmas of the $C^{1+\alpha}$ Pesin theory.  Moreover, the
first approximation property is also valid in the $C^1$ setting with
limit domination by using Liao's shadowing lemma for
quasi-hyperbolic orbit segments\cite{ST}.

The specification property for Axiom A systems ensure the two
approximation properties  in \cite{Sig}. However, the specification
property in a weaker version for non-uniformly hyperbolic systems in
\cite{Hir,LLS,ST} is invalid to the second approximation property,
though it can deduce the first one.
%In the present paper, we develop a new specification property and use it to prove the second result for $C^{1+\alpha}$ non-uniformly hyperbolic
%diffeomorphisms.
More precisely,
%The original technique for Axiom A systems in \cite{Sig} is not suitable for non-uniformly hyperbolic ones.
to achieve the second approximation property, Sigmund\cite{Sig} uses
the specification property infinitely many times to  find the needed
orbit. However, for the nonuniformly hyperbolic case, his idea is not
suitable: the specification property for \textit{finite} orbit segments in the same Pesin block, introduced in \cite{Hir,LLS,ST}, can not be used infinitely many times (even two times),
since we can not determine that the given periodic points and the shadowing
periodic orbits always stay in the required set $\tilde \Lambda$. Therefore, to deal with non-uniformly
hyperbolic case, we disinter a new specification property for \textit{infinite} orbit segments (allowing belonging to different Pesin blocks), inspired
from Katok's Shadowing Lemma,  and use it only once to find the
needed orbit and hence avoid induction. Now we start to introduce
our results precisely.

Throughout this paper, we consider an $f\in \Diff^{1+\alpha}(M)$ and an ergodic hyperbolic
measure $\omega$ for $f$. Let $\Lambda=\cup _{\ell=1}^\infty \Lambda_\ell$ be the Pesin set
associated with $ \omega$. We denote by $\omega
|_{\Lambda_{\ell}}$ the conditional measure of $\omega$ on
$\Lambda_{\ell}.$ Set $\tilde \Lambda_\ell=supp(
\omega|_{\Lambda_{\ell}})$ and $\tilde
\Lambda=\cup_{\ell=1}^\infty\tilde \Lambda_\ell.$ Clearly, $f^{\pm
1}\tilde \Lambda_\ell\subset \tilde \Lambda_{\ell +1},$ and the
sub-bundles $E^s(x),\,\, E^u(x)$ depend continuously on $x\in\tilde
\Lambda_{\ell}.$  Moreover,  $\tilde \Lambda$ is $f-$invariant with
$\omega $-full measure.

We denote by $V_f(\nu)$ the set of accumulation measures of time
averages $$\nu^N=\frac 1{N}{\sum_{j=0}^{N-1}f_\ast^j\nu}.$$ Then $V_f(\nu)$ is a nonempty, closed and connected subset of ${\cal M}_{inv}(M)$. And we denote by
$V_f(x)$ the set of accumulation measures of time averages
$$\nu^N=\frac 1{N}{\sum_{j=0}^{N-1}\delta_{f^jx}},$$ where $\delta_x$ denotes the Dirac measure at $x$. Now we state our main theorems as follows.

\begin{Thm}\label{Thm:maxosc}
For every nonempty connected
set $V \subseteq\{\nu\in{\cal M}_{inv}(M)\,|\,\nu(\tilde{\Lambda})=1\},$ there exists a point $x\in M$
such that $$Closure(V) = V_f(x).$$ Moreover, the set of such $x$ is dense
in $supp(\omega)$, that is, the closure of this set contains $supp(\omega)$.
\end{Thm}
\smallskip

A point $x\in M$ is called to be {\it a generic point} for an $f-$invariant measure $\nu$ if for any $\phi\in C^0(M,\mathbb{R})$, the limit $\lim_{n\to\infty}\frac1{n}\sum^{n-1}_{i=0}\phi(f^ix)$ exists and is equal to $\int\phi d\nu$. As a corollary of Theorem \ref{Thm:specification}, the following holds.

\begin{Cor}\label{Cor:genericpts}
Every $f-$invariant measure supported on $\tilde\Lambda$ has generic points and all generic points form a dense subset in $supp(\omega)$.
\end{Cor}
\smallskip

%\begin{def}\label{defmaxosc}
A point $x\in M$ is said to have {\em maximal
oscillation} if $$V_f(x)\supseteq Closure\{\nu\in{\cal M}_{inv}(M)\,|\,\nu(\tilde{\Lambda})=1\}.$$
%\end{def}
We can deduce from Theorem \ref{Thm:maxosc} that the points having maximal oscillation are dense in $supp(\omega)$. As an extension to Theorem \ref{Thm:maxosc}, we go on to prove that they form a residual subset of $supp(\omega)$.
\smallskip

\begin{Thm}\label{Thmresidual}
The set of points having
maximal oscillation is residual in $supp(\omega)$.
\end{Thm}
\smallskip

\begin{Rem}
For any homeomorphism $f:X\to X$ on a compact metric space preserving an ergodic measure $\omega$, if $(f,\,\omega)$ has specification property(see Theorem \ref{Thm:specification} for more details), analogous arguments and results as in Theorem \ref{Thm:maxosc} and \ref{Thmresidual} are adaptable.\qed
\end{Rem}
\smallskip

A point is called to be an irregular point if there is a continuous function $\phi\in C^0(M,\mathbb{R})$, such that the limit $\lim_{n\to\infty}\frac1{n}\sum^{n-1}_{i=0}\phi(f^ix)$ does not exist. As an application of Theorem \ref{Thmresidual}, we have the below result.
\smallskip

\begin{Thm}\label{Thm:irregular}
If $Closure({\mathcal{M}}_{inv}(\tilde\Lambda))$ is nontrivial(i.e., contains at least one measure different from $\omega$), then the set of all irregular points is residual in $supp(\omega)$.
\end{Thm}
\smallskip

This paper is organized as follows. In section 2, we recall the definition of Pesin set and Katok's shadowing lemma. In section 3, we develop a new specification property and verify that $(f,\,\omega)$ admits this property. In section 4, we use the information on orbit segments to describe that of an invariant measure. In section 5 we use the results in section 3 and 4 to prove Theorem \ref{Thm:maxosc} and then in section 6 we use Theorem \ref{Thm:maxosc} to prove Theorem \ref{Thmresidual} and \ref{Thm:irregular}.

\bigskip

\section{Preliminaries}
\smallskip
We recall the concept of Pesin set and recall some preliminary lemmas in this section.

\smallskip

\subsection{\bf Pesin set (\cite{K3,P1})}
 Given $\lambda,\mu \gg \varepsilon>0$, and for all $k \in \mathbb{Z}^{+}$,
we define $\Lambda_{k}=\Lambda_{k}(\lambda,\mu;\varepsilon)$ to be
all points $x \in M$ for which there is a splitting
$T_{x}M=E_{x}^{s} \oplus E_{x}^{u}$ with invariant property
$D_{x}f^{m}(E_{x}^{s})=E_{f^{m}x}^{s}$ and
$D_{x}f^{m}(E_{x}^{u})=E_{f^{m}x}^{u}$ satisfying:

\begin{enumerate}
\item[(a)] $\|Df^{n}|_{E_{f^{m}x}^{s}}\| \leq e^{\varepsilon
k}e^{-(\lambda-\varepsilon)n}e^{\varepsilon \mid m\mid},~\forall
m\in\mathbb{Z},~n\geq 1$;

\item[(b)] $\|Df^{-n}|_{E_{f^{m}x}^{u}}\| \leq e^{\varepsilon
k}e^{-(\mu-\varepsilon)n}e^{\varepsilon \mid m\mid},~\forall m\in
\mathbb{Z}, ~n\geq 1$;

\item[(c)] $\tan (\angle(E_{f^{m}x}^{s},E_{f^{m}x}^{u})) \geq
e^{-\varepsilon k}e^{-\varepsilon \mid m\mid},~\forall m\in
\mathbb{Z}$.
\end{enumerate}

We set
$\Lambda=\Lambda(\lambda,\mu;\varepsilon)=\bigcup_{k=1}^{+\infty}
\Lambda_{k}$ and call $\Lambda$ a Pesin set.

It is obvious that if $\varepsilon_{1}<\varepsilon_{2}$, then
$\Lambda(\lambda,\,\mu;\,\varepsilon_{1})\subseteq\Lambda(\lambda,\,\mu;\,\varepsilon_{2})$.
\smallskip

According to Oseledec Theorem, $\omega$ has $s \,\,(s\leq
d=dim M)$ nonzero  Lyapunov exponents
$$\lambda_{1}< \cdot\cdot\cdot<\lambda_{r}<0<\lambda_{r+1}<\cdot\cdot\cdot< \lambda_{s}$$
with associated Oseledec splitting
$$T_{x}M=E_{x}^{1}\oplus\cdot\cdot\cdot\oplus E_{x}^{s},\,\,\,\,x\in O(\omega),$$
where we recall that $O(\omega)$ denotes an Oseledec basin of $\omega.$ If
we denote by $\lambda$ the absolute value of the largest negative
Lyapunov exponent $\lambda_r$ and $\mu$ the smallest positive
Lyapunov exponent $\lambda_{r+1}$ and set $E^s_x=
E^1_x\oplus\cdot\cdot\cdot\oplus E^r_x,$ $E^u_x=E^{r+1}_x\oplus
\cdot\cdot\cdot\oplus E^s_x,$ then we get a Pesin set
$\Lambda=\Lambda(\lambda,\mu;\varepsilon)$ for a small
$\varepsilon$. We call it the Pesin set associated with $\omega.$ It
follows(see, for example, Proposition 4.2 in \cite{P1}) that
$\omega(\Lambda\setminus O(\omega))+ \omega(O(\omega)\setminus \Lambda)=0.$

The following statements are elementary:
\begin{enumerate}
\item[(a)] $\Lambda_{1} \subseteq \Lambda_{2} \subseteq
\Lambda_{3}\subseteq \cdot\cdot\cdot$;

\item[(b)] $f(\Lambda_{k}) \subseteq
\Lambda_{k+1},~f^{-1}(\Lambda_{k}) \subseteq \Lambda_{k+1}$;

\item[(c)] $\Lambda_{k}$ is compact for $\forall\,\, k\geq 1$;

\item[(d)] for $\forall\,\, k\geq 1$ £¬the splitting $x\to E_{x}^{u}\oplus
E_{x}^{s}$ depends continuously on  $x\in\Lambda_{k}$.
\end{enumerate}

\bigskip

\subsection {\bf Shadowing lemma}

Let $(\delta_{k})_{k=1}^{+
\infty}$ be a sequence of positive real numbers. Let
$(x_{n})_{n=-\infty}^{+ \infty}$ be a sequence of points in
$\Lambda=\Lambda(\lambda, \mu, \varepsilon)$ for which there exists
a sequence $(s_{n})_{n=-\infty}^{+ \infty}$ of positive integers
satisfying:

\begin{enumerate}
\item[(a)] $x_{n}\in \Lambda _{s_{n}}, ~\forall n\in \mathbb{Z}$;

\item[(b)] $\mid s_{n}-s_{n-1}\mid \leq 1, ~\forall n\in \mathbb{Z}$;

\item[(c)] $d(fx_{n},x_{n+1})\leq \delta_{s_{n}}, ~\forall n\in \mathbb{Z}$;
\end{enumerate}
then we call $(x_{n})_{n=-\infty}^{+ \infty}$ a
$(\delta_{k})_{k=1}^{+ \infty}$  pseudo-orbit. Given $\eta>0$, a
point $x\in M$ is an
 $\eta$-shadowing point for the $(\delta_{k})_{k=1}^{+ \infty}$
 pseudo-orbit if $d(f^{n}x,x_{n+1})\leq \eta \varepsilon_{s_{n}},~ \forall n\in
\mathbb{Z}$, where $\varepsilon_{k}=\varepsilon_{0}e^{-\varepsilon
k}$ and $\varepsilon_0$ is a constant.
\bigskip

\begin{Lem}\label{LemShadow}
(Shadowing lemma \cite{K3,P1}) Let $f:M\rightarrow M$ be a
$C^{1+\alpha}$ diffeomorphism, with a non-empty Pesin set
$\Lambda=\Lambda(\lambda,\mu;\varepsilon)$ and fixed parameters,
$\lambda,\mu\gg\varepsilon>0$. For $\forall \eta >0$ there exists a
sequence $(\delta_{k})_{k=1}^{+ \infty}$ such that for any
$(\delta_{k})_{k=1}^{+ \infty}$ pseudo-orbit there exists a unique
$\eta$-shadowing
point.
\end{Lem}
\bigskip

\section{Specification Property for Non-uniformly hyperbolic systems }

In this section, we develop a new specification property for $C^{1+\alpha}$ non-uniformly hyperbolic systems, which will play crucial role in the proof of Theorem \ref{Thm:maxosc}.
\bigskip

\begin{Thm}\label{Thm:specification} $(f,\,\omega)$ has specification property in the following sense. For any $\eta>0$ there is a sequence of integers
$\{M_{k,\,l}=M_{k,\,l}(\eta)\}_{k,\,l\geq 1}$ satisfying:

Given a sequence of integers $\{k_s\,|\,k_s\geq
1\}_{s\in [a,\,b]\cap\mathbb{Z}}$ for any $-\infty\leq a<b\leq\infty$ and a sequence of orbit segments
$\big{\{}\,\{f^i(x_s)\}_{i=0}^{n_s}\,|\,x_s,\,f^{n_s}x_s\in\,\tilde\Lambda_{k_s},
\,n_s\in\mathbb{N}\,\big{\}}_{s\in [a,\,b]\cap\mathbb{Z}}$, there exist a shadowing point
$z\in M$ and an increasing sequence of integers
$\{c_s\}_{s\in [a-1,\,b]\cap\mathbb{Z}}$ with $0\leq c_{s+1}-c_s-n_{s+1}\leq
M_{k_s,\,k_{s+1}}$ $(s\in [a-1,\,b-1]\cap\mathbb{Z})$ such that
$$d(f^{c_{s-1}+j}z,f^jx_s)<\eta\,\varepsilon_{_{k_{_s}+1}},\,\forall\,j=0,\,1,\,\cdots,n_s-1,\,s\in [a,\,b]\cap\mathbb{Z},$$
where $\varepsilon_k=\varepsilon_0 e^{-\varepsilon k}$ and $\varepsilon_0$ is a constant.

In particular, if $a$ and $b$ are finite integers, the shadowing point $z$ should be periodic with period $\pi=c_b-c_{a-1}$.
\end{Thm}
\smallskip

\begin{Rem}\label{r:explain}The consequence of Theorem 3.4\cite{LLS} or Hirayama's definition for specification property is a particular case of the above theorem. More precisely, they considered {\it finite} orbit segments and asked the beginning and ending points of these segments must be in the {\it same} block $\tilde\Lambda_l.$
\end{Rem}
\smallskip

{\bf Proof of Theorem \ref{Thm:specification}} \,\,

 For $\forall\,\, \eta
>0$, by Lemma \ref{LemShadow} there exists a sequence
$(\delta_{k})_{k=1}^{+ \infty}$ such that for any
$(\delta_{k})_{k=1}^{+ \infty}$ pseudo-orbit there exists a unique
$\eta$-shadowing point.

Let $k_*$ big enough such that $\omega(\tilde\Lambda_k)>0$ for all $k\geq
k_*$. For every $k\geq k_{*}$, take and fix for $\tilde\Lambda_k$
a finite cover $\alpha_k=\{V^k_1,V^k_2,\cdots,V^k_{r_k}\}$ by
nonempty open balls $V^k_i$ in $M$ such that diam$(U^k_i)
<\delta_{k+1}$ and $\omega(U^k_i)>0$ where $U^k_i=V^k_i\cap
\tilde\Lambda_k,\,i=1,\,2,\,\cdots,\,r_k$. Since $\omega$ is $f-$ergodic, by Birkhoff ergodic theorem we have
\begin {equation}\label{Eq:4}
\lim_{n\rightarrow +\infty}\frac 1n
  \sum_{h=0}^{n-1}\omega(f^{-h}(U^\ell_i)\cap U^k_j)=\omega(U^\ell_i)\omega(U^k_j)>0,
\end{equation}
$\forall k,\ell\geq k_\ast$, $\forall 1\leq i\leq r_\ell, 1\leq j\leq r_k$. Then
take
\begin {equation}\label{Eq:11} X^{k,\,\ell}_{i,\,j}=\min\{h\in \mathbb{N}\,\,|\,\, h\geq
 1, \,\,\omega(f^{-h}(U^\ell_i)\cap
  U^k_j)>0\}.
\end{equation}
By  (\ref{Eq:4}), $1\leq X^{k,\,l}_{i,\,j}<+\infty$. Let
$$M_{k,\,\ell}=
\max_{1\leq i\leq r_k,\,1\leq j\leq r_l}X^{k,\,\ell}_{i,\,j}.$$

Now let us consider an increasing sequence of integers $\{k_s\,|\,k_s\geq
k_{*}\}_{s\in\mathbb{Z}}$ and a sequence of orbit segments
$\big{\{}\,\{f^i(x_s)\}_{i=0}^{n_s}\,\,|\,\,x_s,\,f^{n_s}x_s\in\,\tilde\Lambda_{k_s},
\,n_s\in\mathbb{N}\,\big{\}}_{s\in\mathbb{Z}}$. For each $s\in\mathbb{Z}$, we take and fix two integers $s_0$ and $s_1$
%$U^{k_s}_{s_0},\,U^{k_s}_{s_1}\in\alpha_{k_s}$
so that $$x_s\in
U^{k_s}_{s_0},\,f^{n_s}x_s\in U^{k_s}_{s_1},\,s\in\mathbb{Z}.$$ Take
$y_s\in U^{k_s}_{s_1}$ by (\ref{Eq:11}) such that
$f^{X^{k,\,l}_{(s+1)_0,\,s_1}}y_s\in U^{k_{s+1}}_{(s+1)_0}$ for
$s\in\mathbb{Z}$. Thus we get a $(\delta_{k})_{k=1}^{+ \infty}$
pseudo-orbit in M:
$$\cdots \{f^t(x_1)\}_{t=0}^{n_1}\,\cup\,
\{f^t(y_1)\}_{t=0}^{X^{k_1,\,k_2}_{2_0,\,1_1}}\,\cup\,
\{f^t(x_2)\}_{t=0}^{n_2}\,\cup\,\{f^t(y_2)\}_{t=0}^{X^{k_2,\,k_3}_{3_0,\,2_1}}
\,\cup\,\cdots.$$ More precisely, $$x_s,\,\, f^{n_s}(x_s)\in
\tilde\Lambda_{k_s}\subseteq\Lambda_{k_s},\,\,y_s\in
\tilde\Lambda_{k_s}\subseteq\Lambda_{k_s}\,\,\text{and}\,\,f^{X^{k_s,\,k_{s+1}}_{(s+1)_0,\,s_1}}y_s\in
\tilde\Lambda_{k_{s+1}}\subseteq\Lambda_{k_{s+1}},$$ and
$$d(f^{n_s}(x_s),\,y_s)<\delta_{k_s+1},\,\,\,d(f^{X^{k_s,\,k_{s+1}}_{(s+1)_0,\,s_1}}y_s,\,x_{s+1})<\delta_{k_{s+1}+1},\,\forall
s\in\mathbb{Z}.$$ Hence  there exists an $\eta-$shadowing point $z\in
M$ such that
$$d(f^{c_{s-1}+j}z,f^jx_s)<\eta\,\varepsilon_{_{k_{_s}+1}},\,\,\forall\,j=0,\,1,\,\cdots,n_s-1,\,\,\,s\in\mathbb{Z},$$
 where $$c_s=\begin{cases}
 0,&\text{for }s=0\\
 \sum_{j=0}^{s-1}[n_j+X^{k_j,\,k_{j+1}}_{(j+1)_0,j_1}],&\text{for
 }s>0 \\
 -\sum_{j=s}^{-1}[n_j+X^{k_j,\,k_{j+1}}_{(j+1)_0,j_1}],&\text{for $s<0$.}
\end{cases}
 $$
 This ends the proof.
 \hfill $\Box$

\bigskip

\section{Characterizing invariant measures by orbit segments}

It is well-known that for ergodic systems, the time average is the same for almost all initial points and coincides with the space average due to Birkhoff Ergodic Theorem. However, it is not true for general measure-preserving systems (for example, the measure supported on two periodic orbits). Inspired by Ergodic Decomposition Theorem, we prove in the following that the space average can be approximated by the information along finite orbit segments.

Given a finite subset $F\subseteq C^0(M,\,\mathbb{R})$, we denote $$\|F\|=max\{\|\xi\|;\xi\in F\}.$$

\begin{Prop}\label{Prop:approx1}
Suppose $f:X\to X$ is a homeomorphism on a compact metric space and $\nu$ is an $f-$invariant measure. Then for any numbers $\varepsilon>0$,
 any finite subset $F\subseteq C^0(M,\,\mathbb{R})$ and any set $\Delta\subseteq X$ with
  $\nu(\Delta)>{(1+\frac{\varepsilon}{16\|F\|})}^{-1}$, there are a measurable partition $\{R_j\}_{j=1}^b$ of $\Delta$,
  $(b\in\mathbb{Z})$ and
 a positive integer $T$, % and positive rational numbers $\theta_1,\,\theta_2,\,\cdots\,\theta_b$ with $\sum_{j=1}^b\theta_j=1$
  such that for any $x_j\in R_j$ and any integers $T_j\geq T$, $(1\leq j\leq b)$, we have that
$$|\int\xi(x)d\nu-\sum_{j=1}^{b}\theta_j\frac1{T_j}\sum_{h=0}^{T_j-1}\xi(f^h(x_j)) |<\varepsilon, \forall \xi\in F, $$ for any $\theta_j>0$ satisfying $|\theta_j-\frac{\nu(R_j)}{\nu(\Delta)}|<\frac{\varepsilon}{2b\|F\|}$, $1\leq j\leq b.$
\end{Prop}
\smallskip

{\bf Proof}\,\,Let $A=\sup \{ |\xi^*(x)|\,\,\big{|}\,\,x\in Q(f),\,\,\xi\in
F\}.$ Denote by $[a]$ the maximal integer not exceeding $a.$ For
$j=1,\cdots,\, [\frac{32A\|F\|}{\varepsilon\,\|\xi\|}]+1,\,\,\xi\in F,$ set
$Q_j(\xi)=\{ x\in Q(f)|\,\, -A+\frac{(j-1)\varepsilon}{16\|F\|}\,\|\xi\|\leq
\xi^*(x)<-A+\frac{j\varepsilon}{16\|F\|}\|\xi\|\}.$ Let ${\cal B}=\bigvee_{\xi\in
F} \{Q_1(\xi),\,\cdots,\,
Q_{[\frac{32A\|F\|}{\varepsilon\,\|\xi\|}]+1}(\xi)\}$, where
$\alpha\vee\beta=\{A_i\cap B_j\,\,|\,\,A_i\in \alpha,\,\,B_j\in
\beta\}$ for partitions $\alpha=\{A_i\},\,\, \beta=\{B_j\}$. Then
${\cal B}=\{R_j\}_{j=1}^{b} $ is a partition of $Q(f).$
Hence that the positive-measure sets in $\{R_j\cap \Delta\}_{j=1}^b$ form a partition of $\Delta$. For simplicity, we still denote this partition by ${\cal B}=\{R_j\}_{j=1}^{b} $.
 Then
by the definition of ${\cal B}$ and $Q_j(\xi)$ above, we have

\begin{eqnarray*}
& &|\int_\Delta\xi^*(x)d\nu-\sum_{j=1}^{b}\theta_j\xi^*(x_j)|\\
&=&|\int_\Delta\xi^*(x)d\nu-\sum_{j=1}^{b}\frac{\nu(R_j)}{\nu(\Delta)}\xi^*(x_j)|
+|\sum_{j=1}^{b}\frac{\nu(R_j)}{\nu(\Delta)}\xi^*(x_j)-\sum_{j=1}^{b}\theta_j\xi^*(x_j)|\\
&\leq&|\int_\Delta\xi^*(x)d\nu-\sum_{j=1}^{b}\nu(R_j)\xi^*(x_j)|+|\sum_{j=1}^{b}\nu(R_j)\xi^*(x_j)
-\sum_{j=1}^{b}\frac{\nu(R_j)}{\nu(\Delta)}\xi^*(x_j)|\\
&+&|\sum_{j=1}^{b}(\frac{\nu(R_j)}{\nu(\Delta)}-\theta_j)\xi^*(x_j)|\\
&\leq&\sum_{j=1}^{b}\nu(R_j)\max\limits_{y\in R_j}|\xi^*(y)-\xi^*(x_j)|+|\sum_{j=1}^{b}\nu(R_j)\xi^*(x_j)|\cdot(\frac{1}{\nu(\Delta)}-1)+\frac{\varepsilon}{2b\|F\|}\sum_{j=1}^{b}|\xi^*(x_j)|\\
&\leq& \frac{1}{8\|F\|}\cdot\varepsilon\|\xi\|\cdot\sum_{j=1}^{b}\nu(R_j)+A\cdot\sum_{j=1}^{b}\nu(R_j)\cdot\frac{\varepsilon}{16\|F\|}+\frac{\varepsilon A}{2\|F\|}\\
&\leq& \frac{1}{8\|F\|}\cdot\varepsilon\|\xi\|+\frac{\varepsilon A}{16\|F\|}+\frac{\varepsilon A}{2\|F\|}\\
&\leq& \frac{11\varepsilon}{16},\,\, \,\,\xi\in F.
\end{eqnarray*}
For the last inequality, note that $A\leq \|F\|.$

On the other hand, we shall take $T$ large enough such that for all $T_j\geq T,$
$$|\xi^*(x_j)-\frac1{T_j}\sum_{h=0}^{T_j-1}\xi(f^h(x_j))|<\frac{\varepsilon}{16},
\,\,\forall\,j=1,2,\cdots,b,\,\xi\in F.$$
Thus

\begin{eqnarray*}
& &|\int_\Delta\xi^*(x)d\nu
-\sum_{j=1}^{b}\theta_j\frac1{T_j}\sum_{h=0}^{T_j-1}\xi(f^h(x_j)) |\\
&\leq& |\int_\Delta\xi^*(x)d\nu-\sum_{j=1}^{b}\theta_j\xi^*(x_j)|\\
&+&|\sum_{j=1}^{b}\theta_j\xi^*(x_j)-\sum_{j=1}^{b}\theta_j\frac1{T_j}\sum_{h=0}^{T_j-1}\xi(f^h(x_j)) |\\
&<&\frac{11\varepsilon}{16}+|\sum_{j=1}^{b}\theta_j(\xi^*(x_j)-\frac1{T_j}\sum_{h=0}^{T_j-1}\xi(f^h(x_j)) )|\\
&<&\frac{11\varepsilon}{16}+\sum_{j=1}^{b}\theta_j\frac{\varepsilon}{16}\leq \frac{3\varepsilon}{4},
\,\,\xi\in F.
\end{eqnarray*}

Note that $\nu(\Delta)>(1+\frac{\varepsilon}{16\|F\|})^{-1}>1-\frac{\varepsilon}{16\|F\|}$. Hence,%Let $\theta_j=\frac{\nu(R_j)}{\nu(\Delta)},\,1\leq j\leq b.$ Then
\setcounter{equation}{0}
\begin{eqnarray}
& &|\int_{M}\xi(x)d\nu-\sum_{j=1}^{b}\theta_j\frac1{T_j}\sum_{h=0}^{T_j-1}\xi(f^h(x_j)) |\nonumber\\
&=&|\int_{M}\xi^*(x)d\nu-\sum_{j=1}^{b}\theta_j\frac1{T_j}\sum_{h=0}^{T_j-1}\xi(f^h(x_j)) |\nonumber\\
&\leq&|\int_{M}\xi(x)d\nu-\int_{\Delta}\xi(x)d\nu|+|\int_{\Delta}\xi^*(x)d\nu-\sum_{j=1}^{b}\theta_j\frac1{T_j}\sum_{h=0}^{T_j-1}\xi(f^h(x_j)) |\nonumber\\
&\leq&\|\xi\|\cdot(1-\nu(\Delta))+\frac{\varepsilon}{4}\nonumber\\
&\leq&\|\xi\|\frac{\varepsilon}{16\|F\|}+\frac{3\varepsilon}{4}\nonumber\\
&<&\varepsilon, \,\,\xi\in F.\nonumber
\end{eqnarray}
This ends the proof.
\qed
\bigskip

The following lemma is Lemma 3.7 in \cite{LLS}.

\begin{Lem}\label{LemMultiClosing}
Let $f: X\to X$ be a homeomorphism  of a compact metric space
preserving an ergodic measure $\omega.$  Let $\Gamma_j\subset
X$ be measurable sets with $\omega(\Gamma_j)>0$ and for
$x\in\Gamma_j$, let  \
$$S(x, \Gamma_j):=\{ r\in \mathbb{N}\,|\,\, f^{r}x\in \Gamma_j\},$$
$j=1,..., k.$ Take $1>\gamma>0,$ $T\geq 1.$ Then for $\omega-$a.e. $x_j\in \Gamma_j$ there exists $n_j=n_j(x_j)\in S(x_j,
\Gamma_j)$ such that $n_j\geq T$ and
$$0<\frac
{|n_{1}-n_{j}|+...+|n_{j-1}-n_{j}|+|n_{j+1}-n_{j}|+...+|n_{k}-n_{j}|}{\Sigma_{j=1}^k
n_j}<\gamma,$$
 where $j=1,..., k.$
\end{Lem}
\bigskip

\begin{Prop}\label{Prop:approx2}
Suppose $f:X\to X$ is a homeomorphism on a compact metric space and $\nu$ is an $f-$invariant measure. Then for any numbers $\varepsilon>0$,
 any finite subset $F\subseteq C^0(M,\,\mathbb{R})$ and any set $\Delta\subseteq X$ with
  $\nu(\Delta)>{(1+\frac{\varepsilon}{16\|F\|})}^{-1}$, there is a measurable partition $\{R_j\}_{j=1}^b$ of $\Delta$,
  $(b\in\mathbb{Z})$ such that for any positive integer $T$,
  and any recurrence points $x_j\in R_j$, there exist recurrence times $T_j\geq T$, $(1\leq j\leq b)$ satisfying
$$|\int\xi(x)d\nu-\frac1{\sum_{j=1}^{b}{\theta_jT_j}}\sum_{j=1}^{b}\theta_j\sum_{h=0}^{T_j}\xi(f^h(x_j))|<\varepsilon, \forall \xi\in F, $$ for any $\theta_j>0$ satisfying $|\theta_j-\frac{\nu(R_j)}{\nu(\Delta)}|<\frac{\varepsilon}{2b\|F\|}$, $1\leq j\leq b.$
\end{Prop}
\smallskip

{\bf Proof} Take the same partition $\{R_j\}_{j=1}^b$ as in Proposition \ref{Prop:approx1}. By Poincar$\acute{e}$'s Recurrence Lemma, $\nu-$almost every points in $R_j$ are recurrence points. Fix a sequence of recurrence points $x_j\in R_j$ and denote their recurrence time by $T_j$.
By the finiteness of $F$ and Lemma \ref{LemMultiClosing}, we can choose integers $T_j\geq T$ for any $T>0$ such that
$$|\frac{\theta_1(T_1-T_j)+\cdots
+\theta_{j-1}(T_{j-1}-T_j)+\theta_{j+1}(T_{j+1}-T_j)+\cdots+\theta_{b}(T_{b}-T_j)}
{\Sigma_{i=1}^{b}\theta_iT_i}|<\frac{\varepsilon}{16\|F\|}$$ for $j=1,\,2,\,\cdots,\,b.$ Thus

\setcounter{equation}{2}
\begin{eqnarray}
& &|
\Sigma_{j=1}^{b}\theta_j \frac 1{T_j}
\sum_{h=0}^{T_j}\xi(f^h(x_j))
-\frac1{\sum_{i=1}^{b}{\theta_iT_i}}
\sum_{j=1}^{b}\theta_j\sum_{h=0}^{T_j}\xi(f^h(x_j)) |\nonumber\\
&=&| \Sigma_{j=1}^{b}\frac
{\theta_j}{(\theta_{1}+\cdots+\theta_{b})}\frac 1{T_j}
\sum_{h=0}^{T_j}\xi(f^h(x_j)) -\Sigma _{j=1}^{b} \frac
{\theta_{j}}{\theta_{1}T_1+\cdots+\theta_{b}T_{b}}
\sum_{h=0}^{T_j}\xi(f^h(x_j)) |\nonumber\\
&=&|\Sigma_{j=1}^{b}
{\theta_j}\frac{\theta_1(T_1-T_j)+\cdots+\theta_{j-1}(T_{j-1}-T_j)+\theta_{j+1}(T_{j+1}-T_j)
+\cdots+\theta_{b}(T_{b}-T_j)}
{\Sigma_{i=1}^{b}\theta_iT_i}\nonumber\\
&\cdot&\frac 1{T_j}\sum_{h=0}^{T_j}\xi(f^h(x_j)) |\nonumber\\
&\leq& |\Sigma_{j=1}^{b}
{\theta_j}\frac{\theta_1(T_1-T_j)+\cdots+\theta_{j-1}(T_{j-1}-T_j)+\theta_{j+1}(T_{j+1}-T_j)
+\cdots+\theta_{b}(T_{b}-T_j)}
{\Sigma_{i=1}^{b}\theta_iT_i}|\cdot\|\xi\|\nonumber\\
&\leq& \Sigma_{j=1}^{b}
\theta_j|\frac{\theta_1(T_1-T_j)+\cdots+\theta_{j-1}(T_{j-1}-T_j)+\theta_{j+1}(T_{j+1}-T_j)
+\cdots+\theta_{b}(T_{b}-T_j)}
{\Sigma_{i=1}^{b}\theta_iT_i}|\cdot\|\xi\| \nonumber\\
&\leq& \Sigma_{j=1}^{b}
\theta_j\frac{\varepsilon}{16\|F\|}\,\|\xi\|\nonumber\\
&\leq&\frac{\varepsilon}{16},\, \,\,\xi\in F.
\end{eqnarray}
Note that $\nu(\Delta)>(1+\frac{\varepsilon}{16\|F\|})^{-1}>1-\frac{\varepsilon}{16\|F\|}$. Combining with Proposition \ref{Prop:approx1} and inequality $(4.3)$, one deduces that
\setcounter{equation}{2}
\begin{eqnarray} & &|\int_{M}\xi(x)d\nu-
\frac1{\sum_{i=1}^{b}{\theta_iT_i}}\sum_{j=1}^{b}\theta_j\sum_{h=0}^{T_j-1}\xi(f^h(x_j))
|\nonumber\\
&\leq&|\int_{M}\xi(x)d\nu-
\Sigma_{j=1}^{b}\theta_j \frac 1{T_j}
\sum_{h=0}^{T_j-1}\xi(f^h(x_j))
|\nonumber\\
%&\leq&|\int_{M}\xi(x)d\nu-\int_{\Delta}\xi(x)d\nu|+|\int_{\Delta}\xi(x)d\nu-\sum_{j=1}^{b}\theta_j\frac1{T_j}\sum_{h=0}^{T_j-1}\xi(f^h(x_j)) |\nonumber\\
&+&| \Sigma_{j=1}^{b}\theta_j \frac 1{T_j}
\sum_{h=0}^{T_j-1}\xi(f^h(x_j))
-\frac1{\sum_{i=1}^{b}{\theta_iT_i}}\sum_{j=1}^{b}\theta_j\sum_{h=0}^{T_j-1}\xi(f^h(x_j))
|\nonumber\\
&\leq&\varepsilon+\frac{\varepsilon}{16}\nonumber\\
&<&2\varepsilon, \,\,\xi\in F.\nonumber
\end{eqnarray}
Hence we complete the proof.
\qed
\bigskip

\begin{Rem}\label{Rem:FinerPartition}
Through the proof of the previous proposition, one can obtain that the conclusion is suitable for any finer partition of $\{R_j\}_{j=1}^b$.\qed
\end{Rem}

\bigskip

\begin{Prop}\label{Prop:approx3}
Let $\nu$ be an $f-$invariant measure supported on $\tilde\Lambda$. Then for any numbers $\zeta,\,\delta>0$ and
 any finite subset $F\subseteq C^0(M,\,\mathbb{R})$, there are a number $k_\nu\in\mathbb{Z^+}$ and orbit segments $\{z_j,\,fz_j,\,...,\,f^{n_j-1}z_j\}^b_{j=1}$ with $z_j,\,f^{n_j}z_j\in\tilde\Lambda_{k_\nu}$ and $d(f^{n_j}z_j,z_{j+1})<\delta$, $j=1,...,b-1$, satisfying that
$$|\int\xi(x)d\nu-\frac1{\sum_{j=1}^{b}{n_j}}\sum^b_{j=1}\sum_{h=0}^{n_j-1}\xi(f^h(z_j))|<\zeta, \forall \xi\in F.$$
\end{Prop}

\bigskip

{\bf Proof}
Take $k_\nu$ large such that $\nu(\tilde\Lambda_{k_\nu})>{(1+\frac{\zeta}{16\|F\|})}^{-1}$. Applying Proposition \ref{Prop:approx2} with $\Delta=\tilde\Lambda_{k_\nu}$ and $\varepsilon=\zeta$, we obtain a finite partition $\{R_j\}^b_{j=1}$ of $\tilde\Lambda_{k_\nu}$ with $diam R_j<\delta$ and recurrence points $x_j\in R_j$ with large recurrence time $T_j$, $j=1,...,b$ satisfying that $$|\int\xi(x)d\nu-\frac1{\sum_{j=1}^{b}{\theta_jT_j}}\sum_{j=1}^{b}\theta_j\sum_{h=0}^{T_j}\xi(f^h(x_j))|<\varepsilon, \forall \xi\in F, \eqno(4.4)$$ for any $\theta_j>0$ satisfying $|\theta_j-\frac{\nu(R_j)}{\nu(\tilde\Lambda_{k_\nu})}|<\frac{\varepsilon}{2b\|F\|}$, $1\leq j\leq b.$

Recall that  $\omega$ is ergodic and thus for any $1\leq j\leq b$, there is an integer
$X_j\geq 1$ such that $$f^{X_j}R_j\cap R_{j+1}
\neq \emptyset, \,\,1\leq j<b,$$ and $$f^{X_b}R_b\cap R_{1}
\neq \emptyset.$$ Take $y_j\in R_j
$ so that $f^{X_j}y_j\in R_{j+1}
$, $1\leq j<b$ and
$f^{X_b} y_{b}\in R_1$.

For $\zeta$ and $b$, there exists $S\in
\mathbb{N}$ such that
 for any integer $s>S$, we have $0<1/s<\frac{\zeta}{b}$. And then there
exists integers $\bar{s}_1,\,\bar{s}_2,\,\cdots,\,\bar{s}_{b}$ satisfying
$\bar{s}_j/s\leq \frac{\nu(R_j)}{\nu(\tilde\Lambda_{k_\nu})}\leq (\bar{s}_j+1)/s.$ It follows from taking
$s_j=\bar{s}_j\,or\,\bar{s}_j+1\,$ that
$$s=\sum_{j=1}^{b}s_j\,\,\mbox{ and }\,\,\,|\frac{\nu(R_j)}{\nu(\tilde\Lambda_{k_\nu})}-\frac {s_j}{s}|<\frac{\zeta}{2b\|F\|}.$$

Take $T_j$ large enough, such that $\sum^b_{j=1}X_j\ll\sum^b_{j=1}s_jT_j$ and hence it holds that
$$|\frac1{\sum_{j=1}^{b}(s_jT_j+X_j)}\sum_{j=1}^{b}(s_j\sum_{h=0}^{T_j-1}\xi(f^h(x_j))
+\sum_{h=0}^{X_j-1}\xi(f^h(y_j)))
-\frac1{\sum_{j=1}^{b}{s_jT_j}}\sum_{j=1}^{b}s_j\sum_{h=0}^{T_j}\xi(f^h(x_j))|<\zeta.$$
This inequality combing (4.4) with $\theta_j=\frac{s_j}{s}$ implies that

$$|\int\xi(x)d\nu-\frac1{\sum_{j=1}^{b}(s_jT_j+X_j)}\sum_{j=1}^{b}(s_j\sum_{h=0}^{T_j-1}\xi(f^h(x_j))
+\sum_{h=0}^{X_j-1}\xi(f^h(y_j)))|<3\zeta.\eqno(4.5)$$

Let $$z_1=\cdots=z_{s_1}=x_1,\,\,z_{s_1+1}=y_1,$$  $$z_{s_1+2}=\cdots=z_{s_1+s_2+1}=x_2,\,\,z_{s_1+s_2+2}=y_2,$$
$$\cdots\cdots$$ $$z_{\sum^j_{h=1}s_h+j+1}=\cdots=z_{\sum^{j+1}_{h=1}s_h+j}=x_{j+1}$$
$$z_{\sum^{j+1}_{h=1}s_h+j+1}=y_{j+1},$$ $$\cdots\cdots$$ $$z_{\sum^{b-1}_{h=1}s_h+b}=\cdots=z_{\sum^{b}_{h=1}s_h+b-1}=x_{b}$$
$$z_{\sum^{b}_{h=1}s_h+b}=y_{b}.$$
These $\{z_j\}_{j=1}^{\sum^{b}_{h=1}s_h+b}$ are the points we want in the proposition and hence we complete the proof.
\qed
\bigskip

\section{Proof of Theorem \ref{Thm:maxosc}}

In this section, we prove Theorem \ref{Thm:maxosc} by using the specification property developed in section 3 and Proposition \ref{Prop:approx3} in section 4.
\bigskip

{\bf Proof of Theorem \ref{Thm:maxosc}}\,\,
 If
$\{\varphi_j\}_{j=1}^{\infty}$ is a dense subset of $C^0(M,\mathbb{R})$, then
$$\tilde d(\nu,\,m)=\sum_{j=1}^{\infty}\frac{|\int \varphi_j d\nu-\int \varphi_j dm|}{2^j\|\varphi_j\|}$$
is a metric on $\mathcal{M}(M)$ giving the weak$^*$ topology, see
e.g. \cite{Walters}. It is well known that $\mathcal{M}_{inv}(M)$ is
a compact metric subspace of $\mathcal{M}(M)$ in the weak$^*$
topology. For any nonempty closed connected set
$V\subseteq\{\nu\in\mathcal{M}_{inv}(M)|\nu(\tilde\Lambda)=1\}$,
there exists a sequence of closed balls $B_n$ in $\mathcal{M}_{inv}(M)$ with radius
$\zeta_n$ in the metric $\tilde d$ with the weak$^*$ topology
such
that the following holds: \\
(a) $B_n\cap B_{n+1}\cap V \neq \emptyset,$\\
(b) ${\cap_{N=1}^{\infty}\cup_{n\geq N}} B_n=Closure(V),$ \\
(c) $\lim_{n\rightarrow+\infty}\zeta_n=0$.\\
By (a), we take $Y_n\in B_n\cap V$.
% and assume $Y_n$ is an atomic measure supported on the orbit of $x_n\in\tilde\Lambda$.
\bigskip

\begin{Rem}\label{Rem:PerOrb}
In \cite{Sig}, Sigmund assume that $Y_n$ is an atomic measure and thus its information can be characterized by its support(periodic orbit). Hence the remain work is to deal with these periodic orbits by specification property for Axiom A systems. But for our case,  we can not directly take $Y_n$ as an atomic measure( even though this is allowed by \cite{LLS}). The main observation is that the support of these periodic measures may not be contained in $\tilde\Lambda$ and therefore, specification property as in Theorem \ref{Thm:specification} becomes invalid. So we emphasis that $Y_n$ must be in $V$ and thus satisfy $Y_n(\tilde\Lambda)=1$. This allows us to choose pseudo-orbits in $\tilde\Lambda$ whose information can characterize that of $Y_n$ and for which the specification property is valid.\qed
\end{Rem}
\bigskip

Take a finite set
$F_n=\{\varphi_j\}_{j=1}^{n}\subseteq
\{\varphi_j\}_{j=1}^{\infty}$. % with $L_n\uparrow\to+\infty$ such that
%$$\sum_{j=L_n+1}^{\infty}\frac{|\int \varphi^n_j d\nu-\int
%\varphi^n_jdm|}{2^j\|\varphi^n_j\|}\leq
%\sum_{j=L_n+1}^{\infty}\frac1{2^{j-1}}<\zeta_n,\,\,\forall\,\nu,\,m\in\mathcal{M}(M).$$
Let $x_* \in \tilde{\Lambda}$ be given and for any $\delta>0$, let $U_0$ be the open ball of radius
$\delta$ around $x_*$. We have to show that there exists an  $x\in
U_0$ such that $Closure(V)= V_f(x)$.We divide the following proof
into four steps.\\

{\bf Step 1}  An estimation of $Y_n\,\,(n\geq 1).$\\

Let $0<\eta<\frac{\delta}{\varepsilon_0}$ be given and by shadowing lemma we can take and fix $\{\delta_k\}$. Fix $n\in\mathbb{N}$. For $\zeta_n,\,F_n$, by Proposition \ref{Prop:approx3} we choose $k_n=k(Y_n)$ and orbit segments $\{z^n_j,\,fz^n_j,\,...,\,f^{n_j-1}z^n_j\}^b_{j=1}$ with $z^n_j,\,f^{n_j}z^n_j\in\tilde\Lambda_{k_n}$ and $d(f^{n_j}z^n_j,z^n_{j+1})<\delta_{k_{n+1}}$, $j=1,...,b-1$, satisfying that
$$|\int\xi(x)d Y_n-\frac1{\sum_{j=1}^{b}{n_j}}\sum^b_{j=1}\sum_{h=0}^{n_j-1}\xi(f^h(z^n_j))|<\zeta_n, \forall \xi\in F_n.$$
Moreover we can take $k_n<k_{n+1}$ for all $n.$

These segments of orbit segments $\{z^n_j,\,fz^n_j,\,...,\,f^{n_j-1}z^n_j\}^b_{j=1}$ form a `periodic' pseudo-orbit. For simplicity, we can  assume that the `periodic' pseudo-orbit is composed by one orbit segment $\{x_n,\cdots ,f^{p_n-1}(x_n)\}$ with $x_n,\,f^{p_n}(x_n)\in \Lambda_{k_n}$ and $d(x_n,f^{p_n}(x_n))<\delta_{k_{n+1}}.$
Thus, the above inequality can be simplified as
$$|\int\xi(x)d Y_n-\frac1{p_n} \sum_{h=0}^{p_n-1}\xi(f^h(x_n))|<\zeta_n, \forall \xi\in F_n.$$ From this for any $m,$ clearly one has
$$|\int\xi(x)d Y_n-\frac1{m\,p_n} \sum_{h=0}^{m \,p_n-1}\xi(f^{h\,\textrm{mod}p_n}(x_n))|<\zeta_n, \forall \xi\in F_n.\eqno(5.6)$$

\smallskip

{\bf Step 2}  Finding a point $\hat{x}\in U_0$ tracing this pseudo-orbit.

Let $M_n=M_{k_{n-1},k_n}(\eta)$ be numbers defined as in Theorem \ref{Thm:specification}.
Define $$\bar{a}_0=\bar{b}_0=0,$$
$$\bar{a}_1=\bar{b}_0+M_1,\,\,\,\bar{b}_1=\bar{a}_1+2(\bar{a}_1+M_2+p_2)p_1$$
$$\bar{a}_2=\bar{b}_1+M_2,\,\,\,\bar{b}_1=\bar{a}_1+2^2(\bar{a}_2+M_3+p_3)p_2$$
$$\cdots\,\,\,\,\,\,\,\,\,\,\,\,\,\,\,\,\,\,\cdots$$
$$\bar{a}_n=\bar{b}_{n-1}+M_n,\,\,\,\bar{b}_n=\bar{a}_n+2^n(\bar{a}_n+M_{n+1}+p_{n+1})p_n$$
$$\cdots\,\,\,\,\,\,\,\,\,\,\,\,\,\,\,\,\,\,\cdots$$

Using Theorem \ref{Thm:specification} and its proof, we can find a point $\hat{x}\in \Lambda $, $\delta-$close to $x_*$, which $\eta-$shadows the orbit segment $\{x_n,\cdots ,f^{p_n-1}(x_n)\}$ for $m_n=2^n(\bar{a}_n+M_{n+1}+p_{n+1})$ times for all $n$ and runs from $f^{p_n}x_n$ to $x_{n+1}$ with a time lag of no more than $M_{n+1}$.
More precisely, there exist $\{a_n\},\{b_n\}$ with $$a_0=b_0=0,$$ $$b_n=a_n+m_n p_n,\,\,\textmd{and} \,\,a_n-b_{n-1}\leq M_n$$ such that
$$d(f^{j}\hat{x},f^{j\,\textrm{mod}{p_n}}x_n)<\eta \varepsilon_{k_n},\,\,\forall a_n\leq j\leq b_n. \eqno(5.7)$$

\bigskip

\begin{Rem}\label{Rem:diffences}
Note that $a_n\leq\bar{a}_n$, $b_n\leq\bar{b}_n$ and $$b_n-a_n=\bar{b}_n-\bar{a}_n=m_n p_n,\quad a_n-b_{n-1}\leq\bar{a}_n-\bar{b}_{n-1}= M_n.$$ So as $n\rightarrow +\infty,$ $b_n$ and $a_{n+1}$ become much larger than $a_n, M_{n+1},p_n$ and $p_{n+1}.$

The original technique for Axiom A systems in \cite{Sig} is not suitable for non-uniformly hyperbolic ones. Sigmund\cite{Sig} uses the specification property to build inductively a sequence of periodic orbits such that the $n-$th orbit shadows both the $(n-1)-$th orbit and the support of the $n-$th center. In this process the support of the centers and these shadowing periodic orbits are always in the hyperbolic set such that the specification property can be used once by once. Finally, these periodic orbits conjugates to a point $\hat{x} $. However, for the nonuniform hyperbolic case, Sigmund's idea face a difficulty. That is the specification property can not be used once by once, since we can not predetermine the Pesin block in which the shadowing periodic orbits stay. Therefore, to deal with non-uniformly hyperbolic cases, we disinter a new specification property. More precisely, instead of dealing induction, we construct an infinitely many orbit segments, inspired from Katok's Shadowing Lemma. And we apply this property once and for all to find $\hat{x}$ and hence avoid induction.\qed
\end{Rem}
\bigskip

{\bf Step 3} verifying $Closure(V)\subseteq V_f(\hat{x})$.

Let $\nu \in Closure(V)$ be given. By (b) and (c) there exists an increasing
sequence $n_k\uparrow\rightarrow\infty$ such that
$$Y_{n_k}\rightarrow\nu.\eqno(5.8).$$ Let $\xi\in \{\varphi_j\}_{j=1}^\infty=\cup_{n\geq1}F_n$ be given. Then there is an integer $n_\xi>0$ such that for any $n\geq n_\xi$, it holds that $\xi\in F_n$. Denote by $w_\xi(\varepsilon)$ the oscillation
$$\max\{\|\xi(y)-\xi(z)\|\,\,|\,d(y,z)\leq\varepsilon\}$$ and by $\nu_n$ the measure $\delta(\hat{x})^{b_n}$. Thus $$\int \xi d
\nu_n=\frac1{b_n}\sum_{j=0}^{b_n-1}\xi(f^j\hat{x}).\eqno(5.9)$$
 Remark that if $A$ is a finite subset of $\mathbb{N}$, $$|\frac{1}{\# A}\sum_{j\in A}\varphi(f^jx)-\frac{1}{max A+1}\sum^{max A}_{j=0}\varphi(f^jx)|\leq \frac{2(max A+1-\# A)}{\# A}\|\varphi\|\eqno(5.10)$$ for any $x\in M$ and $\varphi\in C^0(M,\mathbb{R})$, where $\# A $ denotes the cardinality of the set $A$.
This inequality (5.10) implies that
$$|\frac1{b_n-a_n}\sum_{j=a_n}^{b_n}\xi(f^j\hat{x})-\frac1{b_n}\sum_{j=0}^{b_n-1}\xi(f^j\hat{x})|
 \leq\frac{2a_n}{b_n-a_n}\|\xi\|,\quad\forall n\geq n_\xi.\eqno(5.11)$$
% Also $$\int \xi d
% Y_n=\frac1{b_n-a_n}\sum_{j=a_n}^{b_n-1}\xi(f^jx_n).$$
On the other hand, combing the inequalities (5.6) and (5.7), one can obtain that $$|\int \xi d
 Y_n-\frac1{b_n-a_n}\sum_{j=a_n}^{b_n}\xi(f^j\hat{x})|\leq \zeta_n+w_\xi(\eta\varepsilon_n),\quad\forall n\geq n_\xi.\eqno(5.12)$$
Note that $$\frac{2a_n}{b_n-a_n}\leq\frac{2\bar{a}_n}{\bar{b}_n-\bar{a}_n}\rightarrow
 0\,\,\, \mbox{ as } n\to \infty, \eqno(5.13)$$ $\zeta_n\to0$ due to Remark \ref{Rem:diffences} and $w_\xi(\eta\varepsilon_n)\rightarrow0$ as $n\to \infty$, it can be deduce by (5.9), (5.11) and (5.12) that $$|\int \xi d
\nu_n-\int \xi d
 Y_n|\leq \zeta_n+w_\xi(\eta\varepsilon_n)+\frac{2a_n}{b_n-a_n}\|\xi\|\rightarrow0,\quad\mbox{ as } n\to\infty.$$ Hence, together with (5.8), it implies that $\nu_{n_k}\rightarrow\nu$ and thus $\nu \in V_f(\hat{x}).$ Therefore, $Closure(V)\subseteq V_f(\hat{x})$.
\bigskip

{\bf Step 4} verifying $V_f(\hat{x})\subseteq Closure(V)$.

Let $\nu\in V_f(\hat{x})$ be given. There exists a sequence $n_k\uparrow\rightarrow\infty$ such that $\nu_{n_k}\to \nu$. Let $\varepsilon>0$ and $\xi\in \{\varphi_j\}_{j=1}^\infty=\cup_{n\geq1}F_n$ be given. For fixed $n_k$, let $i=i(n_k)$ be the largest integer such that $b_{i-1}\leq n_k$. Let $n_k$ (and hence $i$) be so large that $$w_\xi(2^{-i+1})<w_\xi(2^{-i+2})<\frac{\varepsilon}{4}.$$ Let $\alpha=1$ if $b_{i-1}\leq n_k\leq a_i.$ Otherwise, $a_i< n_k\leq b_i$. Write $n_k=a_i+mp_i+l,\,\,0\leq l<p_i$ and define $$\alpha=(b_{i-1}-a_{i-1})(b_{i-1}-a_{i-1}+n_k-a_i-l)^{-1}.$$

Recall $$\int\xi d\nu_{n_k}=\frac{1}{n_k}\sum^{n_k-1}_{j=0}\xi(f^j\hat{x}).$$ Using the inequality (5.10) again, with $A=[a_{i-1}, b_{i-1})\cup [a_{i},\,n_k-l)$, one obtain
\setcounter{equation}{13}
\begin{eqnarray}
& &|\int\xi d\nu_{n_k}-\frac{1}{b_{i-1}-a_{i-1}+n_k-a_i-l}(\sum_{j=a_{i-1}}^{b_{i-1}-1}\xi(f^j\hat{x})
+\sum_{j=a_i}^{n_k-1-l}\xi(f^j\hat{x}))|\nonumber\\
&\leq&2(l+a_{i-1}+a_i-b_{i-1})(b_{i-1}-a_{i-1}+n_k-a_i-l)^{-1}\|\xi\|\nonumber\\
&\leq&2(\frac{p_i}{b_{i-1}-a_{i-1}}+\frac{a_{i-1}}{b_{i-1}-a_{i-1}}
+\frac{M_i}{b_{i-1}-a_{i-1}})\|\xi\|\nonumber\\
&\leq&2(\frac{p_i}{\bar{b}_{i-1}-\bar{a}_{i-1}}+\frac{\bar{a}_{i-1}}{\bar{b}_{i-1}-\bar{a}_{i-1}}
+\frac{M_i}{\bar{b}_{i-1}-\bar{a}_{i-1}})\|\xi\|\nonumber\\
&\leq&\varepsilon \|\xi\|
\end{eqnarray}
provided $n_k$ are large enough due to Remark \ref{Rem:diffences}.
\bigskip

\begin{Rem}
In \cite{Sig}, Sigmund defined $$a_0=b_0=0$$ and $$a_i=b_{i-1}+M_i,\quad\quad b_i=a_i+2^i(a_i+M_{i+1})p_i,\quad i\in\mathbb{N}.$$ It is obvious that these $b_{i-1}$ and $a_{i}$ were chosen independent of $p_{i}$. Here, in our definition(before (5.7)), the choice of $b_{i-1}$ and $a_{i}$  are chosen much larger not only than $a_{i-1}, M_{i},p_{i-1}$ but also than $p_{i}.$ This is one of the important differences to Sigmund's proof. In fact, the assumption of $$n_k-a_i=m p_i$$ in Step 4 in Sigmund's proof is not suitable. The remainder $\ell$ is not greater than $p_i$. However, in his proof, the period $p_i$ may not be small comparing with the lap $b_{i-1}-a_{i-1}$ and hence that $\ell$ is not small enough with respect to $b_{i-1}-a_{i-1}$, which is necessary to the proof as shown in the above inequality (5.14).\qed
\end{Rem}
\bigskip

Then inequality (5.14) implies that
 $$|\int\xi d\nu_{n_k}-[\alpha\frac{1}{b_{i-1}-a_{i-1}}\sum_{j=a_{i-1}}^{b_{i-1}-1}\xi(f^j\hat{x})
+(1-\alpha)\frac{1}{n_k-a_{i}-l}\sum_{j=a_i}^{n_k-1-l}\xi(f^j\hat{x})]|
\leq\varepsilon \|\xi\|.$$
Set $$\rho_{n_k}=\alpha Y_{i-1}+(1-\alpha)Y_i.$$

Using inequality (5.12), one has $$|\int\xi d\nu_{n_k}-\int\xi d\rho_{n_k}|\leq2\varepsilon \|\xi\|$$ for $k$ large enough such that $n_k\gg n_\xi$. Thus $\rho_{n_k}$ has the same limit as $\nu_{n_k}$, that is, $\nu$.

On the other hand, the limit of $\rho_{n_k}$ has to be in $Closure(V)$, since $$\tilde{d}(\rho_{n_k},\,V)\leq \tilde{d}(\rho_{n_k},\,Y_{i})\leq \tilde{d}(Y_{i-1},\,Y_i)\leq 2\zeta_{i-1}+2\zeta_i$$ and $\zeta_i\downarrow0$. Hence, $\nu\in Closure(V)$.

The arbitrariness of $x_*\in\tilde\Lambda$ and $\delta$ implies the density of  $\hat{x}$ in $\tilde\Lambda$. Note that $\tilde\Lambda\subseteq supp(\omega)$ and $\omega(\tilde\Lambda)=1$ and $\omega$ is an ergodic measure. All these conditions ensure that $Closure(\tilde\Lambda)=supp(\omega)$. Hence, it holds that such $\hat{x}$ are dense in $supp(\omega)$. This ends the whole proof.
\qed
\bigskip

\begin{Rem}\, Note that $\m_{inv}(\tilde\Lambda)$ is convex but may not be compact. For better understanding Theorem\ref{Thm:maxosc}, here
we construct a compact connected subset of
$\m_{inv}(\tilde\Lambda)$. Let
$\bar{\varepsilon}=(\varepsilon_1,\varepsilon_2,...,)$ be a (weak)
decreasing sequence of positive real numbers which approach zero.
Let
$$\m_{\bar{\varepsilon}}=\{\nu\in\m_{inv}(f)\,:\,\nu(\tilde{\Lambda}_\ell)\geq1-\varepsilon_\ell,\,\ell=1,2,...\}.$$
Since each $\tilde\Lambda_\ell$ is compact, the map $\nu\to
\nu(\tilde\Lambda_\ell)$ is upper-semicontinuous. Hence,
$\m_{\bar\varepsilon}$ is a closed convex subset of $\m_{inv}(f)$,
the set of all the invariant measures of $f$. This implies
$\m_{\bar\varepsilon}$ is a compact connected subset of
$\m_{inv}(f)$. Since every $\nu\in \m_{\bar\varepsilon}$ satisfying
$\nu(\tilde\Lambda)=1$, we can regard $\m_{\bar\varepsilon}$ as a
subset of $\m_{inv}(\tilde\Lambda)$. Thus, $\m_{\bar\varepsilon}$
must be a compact connected subset of $\m_{inv}(\tilde\Lambda)$.\qed
\end{Rem}
\bigskip

\section{Proof of Theorem \ref{Thmresidual} and  \ref{Thm:irregular}}

In this section, we use Theorem \ref{Thm:maxosc} to prove Theorem \ref{Thmresidual} and then use Theorem \ref{Thmresidual} to prove Theorem \ref{Thm:irregular}.
\bigskip

{\bf Proof of Theorem \ref{Thmresidual}} The proof is not difficult and analogical with the proof
of Proposition 21.18 in \cite{DGS}. Since
$\m(f)$ is compact and convex, we can find open balls $B_n$, $C_n$ in
$\m(f)$ such that
\begin{enumerate}
\item[(a).] $B_n\subset Closure(B_n)\subset C_n$;
\item[(b).] $diamC_n\to 0$;
\item[(c).] $B_n\cap Closure\{\nu\in{\cal M}_{inv}(M)\,|\,\nu(\tilde{\Lambda})=1\}\neq\emptyset$;
\item[(d).] each point of $Closure\{\nu\in{\cal M}_{inv}(M)\,|\,\nu(\tilde{\Lambda})=1\}$ lies in infinitely many
$B_n$.
\end{enumerate}

Put $$P(C_n)=\{x\in M\,|\,V_f(x)\cap
C_n\neq\emptyset\},\quad\forall n\in\mathbb{Z}^+ .$$ It can be
verified that the set of points with maximal oscillation is just
$\cap_{n\geq1} P(C_n)$. Note that
\begin{eqnarray}
P(C_n)&\supseteq&\{x\in M\,|\,\forall N_0\in\mathbb{Z}^+, \,\exists N>N_0 \mbox{ with } \delta(x)^N\in B_n\}\nonumber\\
&=&\cap^\infty_{N_0=1}\cup_{N>N_0}\{x\in\tilde\Lambda\,|\,\delta(x)^N\in B_n\}\nonumber\\
\end{eqnarray}
Since $x\to\delta(x)^N$ is continuous (for fixed $N$), the sets $\cup_{N>N_0}\{x\in M\,|\,\delta(x)^N\in B_n\}$ are open. Since $B_n\cap Closure\{\nu\in{\cal M}_{inv}(M)\,|\,\nu(\tilde{\Lambda})=1\}\neq\emptyset$, these sets are also dense, as shown in Corollary \ref{Cor:genericpts}. Hence $\cap_{n\geq1} P(C_n)$ contains a dense $G_\delta-$set.
\qed

\bigskip

{\bf Proof of Theorem \ref{Thm:irregular}}
Let $x$ be a point having maximal oscillation. By assumption there at least exist two invariant measures $\mu_1\neq\mu_2\in V_f(x).$ So there is continuous function $\phi$ such that $$\int \phi d\mu_1\neq\int \phi d\mu_2.\eqno(6.16)$$ Due to the definition of $V_f(x)$, there are two sequences of integers $n_k,m_k\rightarrow +\infty$ such that $$\delta^{n_k}(x)\rightarrow \mu_1,\,\,\,\,\,\,\delta^{m_k}(x)\rightarrow \mu_2.\eqno(6.17)$$ These imply that $$\lim\frac1{n_k}\sum_{j=0}^{n_k}\phi(f^jx)=\int \phi d\mu_1\mbox{ and }\lim\frac1{m_k}\sum_{j=0}^{m_k}\phi(f^jx)=\int \phi d\mu_2.$$ Combining these equalities with (6.16), we can deduce that $$\lim\frac1{n_k}\sum_{j=0}^{n_k}\phi(f^jx)=\int \phi d\mu_1\neq\int \phi d\mu_2=\lim\frac1{m_k}\sum_{j=0}^{m_k}\phi(f^jx).$$  Thus we have that $\lim\frac1{n}\sum_{j=0}^{n}\phi(f^jx)$
does not exist.
\qed
\bigskip

{\bf Acknowledgement.} The authors thank very much to the whole
seminar of dynamical systems in Peking University.

\bigskip
\end{document}